\documentclass[11pt,twoside]{amsart}
\usepackage[dvips]{graphicx}
\usepackage{amssymb}
\usepackage{latexsym}
\usepackage{xypic}

\newtheorem{theorem}{Theorem}[section]

\newtheorem{proposition}[theorem]{Proposition}

\newtheorem{Thm}{Theorem}

 \newtheorem{question}[Thm]{Question}
\theoremstyle{definition}
\newtheorem{definition}[theorem]{Definition}
\newtheorem{example}[theorem]{Example}

\theoremstyle{remark}

\newenvironment{definition-proposition}{\begin{def-prop} \em}{\end{def-prop}}

\def\P{\mathbb P_{\mathbb C}}
\def\H{\mathbb H}

\numberwithin{equation}{section}

\newcommand{\C}{\mathbb{C}}
\newcommand{\R}{\mathbb{R}}
\newcommand{\s}{\mathbb{S}}

\title{ {An overview of Complex Kleinian Groups} }
\author{Angel Cano \& Jos\'e Seade}
\address{Instituto de Matem\'aticas, Unidad Cuernavaca, \newline
 Universidad  Nacional Aut\'onoma de M\'exico, \newline
 Lomas de Chamilpa, Cuernavaca, M\'exico.  }

\email{angel@matcuer.unam.mx, jseade@matcuer.unam.mx}

\thanks{Research  supported by grants from CNPq   }
\keywords{Equicontinuity, quasi-projective transformations, complex hyperbolic space, discrete groups, limit set}

\subjclass{Primary: 3nQ45, 37F45; Secondary nnE40, 57R30}

\dedicatory{Dedicated to Vale, on the occasion of his 65th birthday.}

\begin{document}

\maketitle

\begin{abstract}
Classical Kleinian groups are discrete subgroups of $PSL(2,\C)$ acting on the complex projective line $\P^1$, which actually coincides with the Riemann sphere, with non-empty region of discontinuity. These  can also be regarded as the monodromy groups of certain differential 
equations. These groups have played a major role in many aspects of mathematics for decades, and also in physics.
It is thus natural to study discrete subgroups of the projective group $PSL(n,\C)$, $ n > 2$.  Surprisingly, this is a branch of mathematics which is in its childhood, and in this article we give an overview of it.
\end{abstract}

\section*{Introduction} 

The study of discontinuous group actions on the complex line   $\Bbb{P}^1_{\Bbb{C}}$ was begun in the 1880's by Schottky,
Poincar\'e and Klein, see \cite{poincare}.   This study  was initially  motivated  by the fact that these groups appeared as  the  
monodromy groups of  Riccati differential equations,  equations  which had shown to be  useful in applications in several problems of mathematical physics, quantum mechanics, and fluid dynamics (see  \cite{levi}). 

The term ``Kleinian group" was coined by Poincar\'e, meaning by this a discrete group of M\"obius transformations acting on the extended complex plane $\widehat \C$ so that no orbit is dense. 
 In fact it can be shown that each Kleinian group can be realized as the monodromy group of a Ricatti differential 
equation. There have been significant developments in the theory coming from various sources, as for instance 
  the theory of quasi-conformal mappings
by Ahlfors and Bers see \cite{crash}, the combination theorems  by  Maskit,   which are  direct descents from previous work by Klein and Koebe see \cite{maskit}. This theory has also strong relations  with the geometrization of 3-manifolds, as remarked by Thurston, see 
\cite{thurston}. There are also remarkable similarities between Kleinian groups and the iteration theory of rational maps, as explained through the Sullivan-MacMullen dictionary.

Around  1980, Yoshida \cite{Yo}  showed that the monodromy groups of the  so-called  Orbifold 
Uniformazing  Differential equations (a generalization of the Ricatti differential equation, now in the partial differential 
equations setting), are  discrete groups  of projective transformations in several complex variables. As in the one dimensional case, a natural 
task is to study  the geometry and  dynamics  of such groups.

This is of course also important from various other points of view. For instance at the end of the 19th Century, E. Picard started  the study of complex hyperbolic geometry,  then continued by other mathematicians including Mostow and Deligne. The projective  Lorenz group $PU(n,1)$ is the group of holomorphic isometries of the complex hyperbolic n-space $\H^n_\C$ , and this is naturally a subgroup of the projective group $PSL(n+1,\C)$. Complex hyperbolic geometry is nowadays a very active field of mathematics, with interesting new articles being published continuously. And by looking at discrete groups of isometries of $\H^n_\C$ one is naturally studying a specially interesting type of discrete subgroups of $PSL(n+1,\C)$. But there are many others, as for instance groups of affine transformations in $\C^n$, and as we will see, there are also Schottky groups, groups appearing through twistor theory, groups appearing through the classical Pappus' theorem, and many others.

The study of discrete subgroups of $PSL(n+1,\C)$ acting on the projective space $\P^n$ with a non-empty region of discontinuity is also important for the uniformisation problem of complex manifolds, since this is a natural way for constructing compact complex manifolds equipped with a projective structure. This was the motivation for  Nori's work in  \cite{No} as well as for several articles by Kato (see \cite{Ka1,Ka2,Ka3,Ka4}).

From the viewpoint of dynamics, possibly the first results for $n>1$ were by Chen and Greenberg in \cite{ChG} for discrete subgroups of $PU(n,1)$. In  \cite{SV1,SV2,SV3} 
Seade and Verjovsky  introduced several ways of constructing discrete subgroups of   $PSL(n+1,\C)$ with very rich geometry and dynamics. Three of these constructions are specially interesting. The first is the {\it suspension construction}, later extended in \cite{Nav1, cano2}. This provides means for constructing subgroups of $PSL(n+1,\C)$ whose geometry and dynamics is governed by a discrete subgroup of $PSL(n+1,\C)$. This is specially interesting when $n=1$ since there is a vast  knowledge and literature about discrete subgroups of $PSL(2,\C)$.

The second main construction in \cite{SV1,SV2} uses twistor theory. This theory is one of the jewels of mathematics and physics in the late 20th Century. Its mathematical foundations were laid down by Penrose and several other mathematicians, including Atiyah, Singer and Hitchin. In a simple way, we can say that this theory associates to each $2n$-dimensional sphere $\s^{2n}$ (we restrict to this setting for simplicity) its twistor space ${\frak Z}(\s^{2n})$, which is a complex projective manifold which fibres over $\s^{2n}$ with fibre at each $x \in \s^{2n}$ the space of all complex structure on the tangent space $T_x \s^{2n}$ which are compatible with the metric and the orientation. In fact the fibre turns out to be the twistor space ${\frak Z}(\s^{2n-2})$. For $n=1$ the twistor space is $\s^2 \cong \P^1$  itself. For $n=2$ the space ${\frak Z}(\s^{4})$ is $\P^3$. Then ${\frak Z}(\s^6)$ is a regular  quadric of complex dimension 6 in $\P^7$, etc.

An starting point for Penrose's twistor program is the fact that there is a rich interplay between the conformal geometry of the spheres and the complex (holomorphic) geometry of their twistor spaces. Then, in \cite{SV1,SV2} the authors show that there is also a rich interplay between conformal dynamics on $\s^4$ and holomorphic dynamics on its twistor space $\P^3$.

The third construction in \cite{SV1,SV3} extends to higher dimensions a classical construction of Schottky groups, thus obtaining a very interesting class of discrete subgroups of automorphisms of $PSL(2n+2,\C)$.

Yet, these constructions essentially work only in odd dimensional projective spaces. So it is natural to ask what happens on even dimensions, particularly in $\P^2$. Surprisingly, this is a branch of mathematics which had essentially passed unnoticed until a few years ago, except for the case of groups appearing in complex hyperbolic geometry. The study  of discrete subgroups of $PSL(3,\C)$ in general somehow started just a few years ago through a 
 series of articles by J. P. Navarrete, A. Cano and W. Barrera (see the bibliography below).

In this expository article we start by reviewing the classical case of Kleinian groups acting on the Riemann sphere $\s^2$, which is isomorphic to the projective line $\P^1$. We then motivate the discussion in higher dimensions by looking at an explicit example in dimension two, where we see that unlike the classical case, in higher dimensions
there is not a ``well-defined" notion of the limit set of a discrete subgroup. There are actually several possible definitions, each with its own interest and characteristics. We then introduce  complex kleinian groups in higher dimensions, a concept introduced in \cite{SV1} which simply means a discrete subgroup of $PSL(n+1,\C)$ which acts properly discontinuously on some non-empty open invariant subset of $\P^n$. An important point here is the definition of the Kulkarni limit set and the Kulkarni region of discontinuity of discrete subgroups of $PSL(n+1,\C)$, which play key-roles in this theory.

Sections 4 and 5 focus on the case $n=2$ and they are mostly based on \cite{cano2}. 
In Section 4 we describe briefly all the known types of complex Kleinian subgroups of $PSL(3,\C)$. In Section 5 we describe with more detail the complex Kleinian subgroups $\Gamma$ of $PSL(3,\C)$ whose action on $\P^2$ has an open invariant region $\Omega$ where the action is properly discontinuous and the quotient $\Omega/\Gamma$ is compact; we call these quasi-cocompact groups. 

Section 6 is based on \cite{BCN, BCN2} and it
springs from the results in \cite{Nav1, Nav2, cano2} studying the structure of ``the limit set" of discrete subgroups  of 
$PSL(3,\C)$. One has that  the Kulkarni limit set always contains at least one projective line, and if it has ``sufficiently many" lines, then its complement, which is the Kulkarni region of discontinuity is the largest open invariant set where the group acts properly discontinuously, and it coincides with the equicontinuity region of the group. Then we explore the cases where ``the limit set has few lines".

Finally, in Section 7 we discuss briefly the situation in higher dimensions. We sketch the constructions of Seade and Verjovsky of Schottky groups and complex Kleinian groups obtained via twistor theory.

\section{Classical Kleinian groups: the one-dimensional case}

  The study of holomorphic  dynamical systems  in one complex variable  relies on  deep connections existing 
between  conformality and holomorphy in the Riemann sphere. In fact, if we consider maps between open subsets of the Riemann sphere, we find that 
every holomorphic function  is conformal and {\it viceversa, i.e.,} every orientation preserving conformal map is holomorphic.
So we start this discussion by considering  conformal groups of  transformations in the Rieman sphere $\s^2$.

We think of  $\s^2$ as being the extended complex line $\widehat \C:= \C \cup \infty$, and 
recall that a conformal map in 
$\s^2$ is a diffeomorphism that preserves angles. 
The simplest conformal maps  are the   reflections  on lines and circles. Recall that such transformations are 
defined as follows:

\begin{definition}  
A {\it reflection} in $\s^2$ on a line $\ell$ is the transformation that associates to each point $(x,y)$ the 
unique point $(x',y')$ such that the segment 
 $\overline{(x,y),(x',y')}$, joining these two points, meets orthogonally the line $\ell$ at a point $P$ which is the middle 
point of this segment. \end{definition}


It is clear that every reflection is an involution, {\it i.e.}, it has order 2. So iterating one of these maps is not interesting.
Let us see what happens when we consider a couple of reflections:

\begin{theorem}
The composition of two reflections on lines which intersect at a point $p$ in an angle $\theta$ is    a rotation 
around $p$ with an angle $2\theta$.
\end{theorem}

\begin{theorem}
The composition of two reflections on  parallel lines is the translation  determined by a vector whose direction is 
orthogonal  to the lines (the sense depends on the order that  we have made the composition) and whose norm is 
$2d$, where $d$ represents the distance between the lines.
\end{theorem}

In the same  line of thought, let us consider now 3 non-concurrent lines 
 $\ell_1 $, $\ell_2,\ell_3$ in $\C$, which determine a  triangle $T$. We assume the angles of $T$  are of the form  
$\pi/p, \,\pi/q,\, \pi/r$. Notice one necessarily has $p, q, r \ge 2$ and $\pi/p + \pi/q + \pi/r = 1$.

Let us denote by  $\Sigma_{p,q,r}$ the group generated by such inversions; this is called 
 a {\it triangle group}. We observe that the orbit of $T$ under $\Sigma_{p,q,r}$ consists of symmetric triangles that cover the 
whole plane $\C$.  
More generally, given any natural numbers $p,q,r \ge 2$ and angles $\pi/p, \,\pi/q,\, \pi/r$, there is 
  exactly one of the classical two-dimensional geometries (Euclidean, spherical, or hyperbolic) which admits a triangle 
$T$, bounded by geodesics,  with angles $(\pi/p, \pi/q, \pi/r)$ (see \cite{ball, cox1, cox2, krai}). 

When the geometry is spherical, the geodesics are segments of circles in $\s^2$ of maximal length, so they are equators. Each such circle corresponds to the intersection of the unit sphere in $\R^3$ with a plane thorough the origin in $\R^3$. One thus has a corresponding reflection on $\R^3$ which determines by restriction a ``reflection" in the 2-sphere.

When the geometry is hyperbolic, the geodesics are segments of circles or lines in $\C$ that meet orthogonally the boundary of the unit disc, which serves as model for hyperbolic geometry. Or equivalently, we can take as a model for complex geometry the upper hemisphere of the unit sphere in $\R^3$. In this case the geodesics are segments of circles which meet orthogonally the equator $\{(x,y,0) \, | \, x^2 + y^2 = 1 \}$.

In all cases the 
corresponding ``plane" can be tiled by the copies of the triangle, obtained by ``reflections" on the three edges of $T$. 
To make this precise we need to extend the type of transformations we consider:


\begin{definition}
A {\it  reflection (or inversion)} in the Riemann sphere  with respect to a circle of center $c$ and radius $r$ is the 
transformation  that carries  each point $p\in \widehat{\Bbb{C}}\setminus \{c,\infty\}$ into the point $q$ such that  
$p,q$ and $c$ are collinear and $\vert c-p\vert \vert c -q \vert=r^2$. The center is sent to   $\infty$ and $\infty$ goes to the 
center of the circle.
\end{definition}


Of course these concepts extend easily to higher dimensions. So, the 3-sphere can be regarded as being $\R^3 \cup \infty$, and we have reflections on 2-planes as well as on 2-spheres. The extension to $\s^n \cong \R^n \cup \infty$ is obvious.

Even though    these transformations are very simple, one has: 

\begin{theorem}
The group $Conf( \Bbb{S}^2)$ of all  conformal automorphisms of the Riemann sphere coincides  with the group 
generated by the reflections on all lines and  circles in the Riemann sphere. Moreover one has: 
\[
Conf( \Bbb{S}^2)=
\langle \{\overline{\mbox{z }}\}\cup M\ddot{o}b(\widehat{\Bbb{C}}) \rangle
\]
where $M\ddot{o}b(\widehat{\Bbb{C}})$ is the M\"obius group consisting of all transformations of the form:
\[
h(z)=\frac{az+b}{cz+d} \;,
\]
with $a,b,c,d$  complex numbers such that $ad-bc=1$.
 \end{theorem}

In other words, every M\"obius transformation is a conformal map, in fact a composition of inversions, and 
the conformal group is generated by the M\"obius group and the complex inversion $z \mapsto \bar z$.


After making this short digression, let us consider again the triangle groups $\Sigma_{p,q,r}$.   These were defined above  when the corresponding triangle was Euclidian. When the triangle is spherical, one has the corresponding group of reflections in $\R^3$, which actually corresponds to a group of reflections of the 2-sphere $\s^2$. Notice that this is the same as considering the corresponding  inversions on $\s^2$.

When the triangle is hyperbolic, we have a corresponding triangle group generated by inversions on the three circles that determine the sides of the triangle.

In all cases
the sum 
of the angles of the triangle determines the type of  geometry one has: this is sustained   by the Gauss-Bonnet theorem. It is 
Euclidean if the  sum of the angles is exactly $\pi$, spherical if it exceeds $\pi$ and hyperbolic if it is strictly smaller 
than $\pi$. Up to permutation of  the numbers $p,q,r$,  which are all $\ge 2$, there are the following possibilities:

\medbreak

\noindent
{\bf a) The Euclidean case.} This happens when  $p^{-1}+q^{-1}+r^{-1}=1$ and 
up to permutations the only possible triples are  $(2,3,6), (2,4,4), (3,3,3)$. The corresponding triangle groups are 


\medbreak

\noindent
{\bf b) The spherical case.} We now have $p^{-1}+q^{-1}+r^{-1}>1$ and 
the triangle group is the finite  group of  symmetries of a tiling of a unit sphere by spherical triangles. The possible triples  are  $(2,3,3), (2,3,4), (2,3,5),$ and $(2,2,r)$, for $ r\ge 2$. Spherical triangle groups can be identified with 
the symmetry groups of regular polyhedra in the three-dimensional Euclidean space: the group induced by the triple 
$(2,3,3)$ corresponds to the tetrahedron, the group $\Sigma_{2,3,4}$ corresponds to both the cube and the 
octahedron (which have the same symmetries group), the triple $(2,3,5)$ corresponds  to both the dodecahedron and 
the icosahedron. 

The groups  $(2,2,n), n\ge 2$, are the symmetries groups of the family of dihedrons, which are 
degenerate solids formed by two identical regular $n$-gons joined together along their boundary, which is a regular $n$-sided polyhedron. We refer for instance to Chapter 2 in \cite{Seade} for a careful description of the relation between triangle groups and the symmetry groups of regular polyhedra.

The spherical tilings corresponding to 
the previous triples  are depicted below:


\medbreak

\noindent
{\bf c) The hyperbolic case.}  We have  $p^{-1}+q^{-1}+r^{-1}<1$.
The triangle group is the infinite  group of  symmetries of  a tiling of the hyperbolic plane by hyperbolic triangles. There are 
infinitely many such groups. Below we depict the  tilings associated with some small values.


Another remarkable fact is that these tilings  appear  naturally  when we consider the  image of the Schwartz  $s$-function of certain hypergeometric differential equations. Moreover, in this case the group formed by the words of even length naturally 
arises as the monodromy group of the respective differential equation, see \cite{love}.

\medbreak

It is natural to ask
what happens if we consider more than three lines?

\medbreak

In the Euclidean case, there is not much more since one has 
the following result, which  is usually attributed  to Leonardo Da Vinci:

\begin{theorem}
The only reflection groups in  Euclidean geometry
which give  rise to   tilings in the plane are:
\begin{enumerate}
 \item The triangular groups: 
$\Sigma_{(3,3,3)}$, $\Sigma_{(2,3,6)}$, $\Sigma_{(2,4,4)}$; 

\item
 The quadrangular group $\Sigma_{(2,2,2,2)}$. 
\end{enumerate}
 
\end{theorem}

In the case of spherical geometry we are talking about discrete subgroups of $SO(3)$, and the only such groups are the triangle groups, described above, and the cyclic groups of finite order. However, in the hyperbolic case the possibilities are infinite. In particular,  given arbitrary integers $p_1,...,p_r \,\ge 2$ such that
$$\pi/{p_1} + \pi/{p_2}  + .... + \pi/{p_r}  < (r-2)\, \pi \,,$$
we have  polyhedrons in the hyperbolic plane, bounded by geodesics, with these interior angles. Each such polyhedron determines a corresponding group of reflections, and a tiling of the hyperbolic plane. Moreover, these are just a special type of discrete subgroups of M\"obius transformations, there are many others.

As we  said before, historically Kleinian groups appeared as a way to studying  in a global way the  solutions of Ricatti 
differential equations. Following this motivation, from now on we are  going to center our interest in  subgroups of $M
\ddot{o}b(\widehat{\Bbb{C}})$. These are  
conformal maps  of the Riemann sphere which preserve the orientation.  Moreover, since the  solutions of the Ricatti differential equations  can be interpreted  as 
Riemann surfaces where the monodromy group corresponds to the monodromy  of the Riemann surface, we are 
motivated  to studying  the set where such groups act properly discontinuously. This is the discontinuity region of the group, and the space of orbits in this region is the corresponding Riemann surface.

 We consider first  the cyclic groups of $M\ddot{o}b(\Bbb{C})$, {\it i.e.}, the groups generated by a single M\"obius 
transformation. The following well-known result describes the dynamics in these cases:
 
\begin{theorem}
Each M\"obius transformation is conjugate to one of the following transformations:

\begin{enumerate}
 \item A  translation  $h(z)=z+1$; in such case the transformation is called  parabolic.
 \item A homothety $h(z)=az$ where $\vert a\vert \neq 1 $;  such  transformations are  called  loxodromic.
\item A rotation $h(z)=az$ where $\vert a \vert =1$; these transformations are known as elliptic.
\end{enumerate}
 \end{theorem}

  It is quite simple to observe that 
 an infinite and   discrete  group  $\Gamma$ of M\"obius transformations  cannot act properly discontinuously on all of $
\Bbb{P}^1_{\Bbb {C}}$, since every infinite sequence in a compact set must have accumulation points. 
 The  {\it limit set} of such a group $\Gamma$, in symbols   $\Lambda(\Gamma)$, is by definition the set of accumulation points
of the $\Gamma$-orbits of all points. Its complement  $\Omega(\Gamma)=\Bbb{P}^1 _{\Bbb{C}}\setminus 
\Lambda(\Gamma)$
is the largest open set where $\Gamma$ acts properly discontinuously.  This  is also the largest open set where  $
\Gamma$  forms a normal family. This important  property   enables us  to establish strong similarities  
with the dynamics of rational functions in  the Riemann sphere, as shown through 
the so-called Sullivan Dictionary (see for instance \cite{Sul}).

\begin{definition}  
A {\it Kleinian group} is a group $\Gamma$ of M\"obius transformations with $\Omega(\Gamma)\neq \emptyset$.  
\end{definition}
 
 In other words, the limit set of a Kleinian group cannot be the whole Riemann sphere.
 
 Notice that 
 this implies that the group must be discrete. 
 
 It is worth saying that this is the classical definition of a Kleinian group. In 
the modern literature, the term ``Kleinian group" often refers to all discrete subgroups of M\"obius transformations. 

It is not hard to show that given a Kleinian group, its limit set  $\Lambda(\Gamma)$ is a closet invariant set which 
either has finite cardinality, consisting of 1 or 2
points, or else it is a nowhere dense perfect
set. In the first case the group is said to be elementary. In the second case the group is nonelementary, and in these 
cases one can show that the limit set is the set of accumulation points of every orbit ({\it i.e.}, independently of the 
choice of orbit).

\medbreak

Let us  proceed to build up a class of Kleinian groups with interesting dynamics called {\it Schottky groups}.

\begin{example} 
Consider  $2g$ , $g\geq 2$,  disjoint closed euclidean discs     $R_1,\ldots,R_g$, $S_1,\ldots,S_g$  in
$\widehat{\Bbb{C}}$,  and   M\"obius transformations  $\gamma_1,\ldots \gamma_g$ such that  
  $\gamma_j(R_j)=\widehat {\Bbb{C}}\setminus \overline{S_j}$. Set $\Gamma=<\gamma_1,\ldots,\gamma_g>$ be the group generated by these maps. Then $
\Gamma$ is a free group with $g$ generators, $\Lambda(\Gamma)$ is a Cantor set and $\Omega(\Gamma)/\Gamma $ is a  
sphere with $g$ handles  attached.
\end{example}


Schottky groups play a significant role in the theory of Kleinian  groups by various reasons. One of these is K\"obe's retrosection theorem, which states 
that every closed oriented 2-manifold can be  uniformized by a  Schottky group, see \cite{maskit}. 

When we allow that  the  discs in 
our definition above  of the Schottky groups touch each other  ``a little", we get the so-called Kissing Schottky groups, beautifully described in \cite{indra}. When we  allow  that the 
discs  overlap, it becomes hard to give conditions to guarantee that  the group is discrete. In this case an important tool is  Maskit's  combination 
theorem,  which allows us to create new dynamics from simpler ones. Another  tool, introduced by Ahlfors and Bers,  which enables us to create new 
examples is by means of the study of  quasi-conformal deformations of a group.


M\"obius transformations can also be regarded as isometries of the real hyperbolic 3-space $\H^3$ 
via Poincar\'e's extension, that we now explain. This provides an important link between Kleinian groups and 3-manifolds theory. In fact Poincar\'e's extension can be done also in higher dimensions, though the way we do it here, which is particularly nice and simple, only works in the dimensions we now envisage.

To explain this, Consider first the quaternionic line (space), which we can think of as being $\R^4$ together with multiplication by the quaternions $i, j, k$. 
Consider   the half-space  $\Bbb{H}^ 3= \{x +i y+ tj \,\vert \,  t > 0\}$ in the quaternionic line. This is conformally equivalent to a 3-ball. By 
identifying the boundary of $\Bbb{H}^ 3$ with $\Bbb{P}^ 1_{\Bbb{C}}$ we can extend the action of $M\ddot{o}b(\hat{\Bbb{C}})$ 
on $ \Bbb{P}^1_{\Bbb{C}} $ to an action on $\overline{\Bbb{H}^ {3}}$  by:
\[
h(w) =(aw+b)(cw+d)^ {-1}.
\]
where $w\in \Bbb{H}^ 3$ and $a,b,c,d$ are complex numbers satisfying $ad-bc= 1$, the multiplication here is the usual multiplication by quaternions.
This  is usually known as the Poincar\'e extension.

It is possible to show the existence of a metric  compatible with the topology of $\Bbb{H}^ 3$ for which each orientation preserving isometry 
is exactly the  Poincar\'e extension of a M\"obius transformation. This is the hyperbolic metric. For any given discrete group of 
$PSL(2,\Bbb{C})$, one has that its action on $\Bbb{H}^ 3$, being by isometries, is properly discontinuous, and we
may consider the quotient space $\Bbb{H}^ {3}/\Gamma$. This is an orbifold and if $\Gamma$ is Kleinian, then  there 
is actually  an associated 3-orbifold  with boundary
  \[
\mathcal{M} (\Gamma)= (\Bbb{H}^ {3}\cup \Omega(\Gamma))/\Gamma.  
  \]
 The boundary $\partial (\mathcal{M}(\Gamma))$ carries a complex structure,  
so it as a disjoint union of Riemann surfaces with an  orbifold structure,
\[
\Omega(\Gamma)/\Gamma = S_1 \cup  S_2 \cup \ldots 
\]
The connected components of $\Omega(\Gamma)$ fall into  the corresponding $\Gamma$-conjugacy classes 
$\{\Omega_i\}$  
covering the various $S_i$.  If we denote by $Isot(\Omega,\Gamma)$ the stabilizer in $\Gamma$ of a given
component $\Omega_i$, then of course $\Omega_i/Isot(\Omega,\Gamma) = S_i$. Such a subgroup is said to be a 
component
subgroup of $\Gamma$. The  celebrated Finiteness Theorem due to Ahlfors asserts that if the group is finitely generated, then $\Omega(\Gamma)$ has a finite number of components and each  one is a Riemann surface of finite volume. Moreover   a key
result due to Bers, known as the Simultaneous Uniformisation Theorem,  states that any finite set of Riemann surfaces with  finite volume can be represented
in this way. \\

\begin{example}
\begin{enumerate}
\item 
If $\Gamma = \langle T \rangle $ with $T$ loxodromic, then  $\mathcal{M}(\Gamma)$ is a solid torus.

\item If $\Gamma = \langle T \rangle $  with $T$
parabolic, then $\mathcal{M}(\Gamma)$ is homeomorphic to \\ $\{0 < \vert z\vert  \leq 1 \} \times  (0,1)$  with boundary
a twice-punctured sphere.

\item If $\Gamma=\langle h(z)=z+w_1, g(z)=z+w_2 \rangle$, then $\mathcal{M} (\Gamma)$ is homeomorphic to
 $\{0 < \vert z\vert  \leq 1 \} \times  \Bbb{S}^ 1$ and $\partial \mathcal{M}(\Gamma)$  is a torus.
 
\item   If $\Gamma$ is a torsion free Fuchsian group whose  limit set is 
$\Bbb{R}\cup \infty $, then $\mathcal{M}(\Gamma)$ is the product  of $\Omega(\Gamma)/\Gamma$ with a closed interval,.
\item  Given a Schottky group, then  $\mathcal{M}(\Gamma)$  is a handle-body with boundary a single compact surface of genus $n$.
\end{enumerate}
\end{example}

Notice that the group of M\"obius transformations can be identified with the projective group 
$$PSL(2,\C) \cong SL(2,\C) / \pm Id \, ,$$
via the map that carries each matrix $
\left (
\begin{array}{ll}
a & b \\
 c & d
\end{array}
\right )
$
into the M\"obius transformation $z \mapsto \frac{az +b}{cz + d}$. 
In other words, we may let $SL(2,\C) $ act linearly on $\C^2$ in the usual way. This action carries lines into lines and therefore determines an automorphism of the complex projective line $\P^1$, which is isomorphic to the Riemann sphere $\s^2$. The induced action on $\s^2$ determines  the corresponding M\"obius transformations, which is as above.

When we are interested  in the generalization of Kleinian groups to higher dimensions, we can see that there are at least  two natural ways to proceed: we may consider either conformal transformations of the $m$-sphere, or else studying holomorphic automorphisms of the complex projective space $\P^n$. When $m=2$ and $n=1$ the two theories coincide because every holomorphic map in one complex variable is conformal and conversely, every orientation preserving conformal map in $\widehat \C$ is holomorphic.
But this does not happen in higher dimensions: neither holomorphic maps  are necessarily conformal, nor conformal maps are necessarily holomorphic.

The study 
 of conformal transformations in higher dimensional 
spheres is most classical. This corresponds to isometries of higher dimensional real hyperbolic spaces, and this is a subject in which there is a vast literature and remarkable works by many authors. A nice overview  of the theory of Kleinian groups  before 1960 can be found in the treatise of
Fatou, see \cite{fatou}, and a good account of the modern study of  conformal Kleinian groups  can be 
found in \cite{kapovich}.

The study of groups of automorphisms of  complex projective spaces includes the theory of lattices in complex hyperbolic geometry, which is a subject that goes back to the work of E. Picard and others. Yet, this is a theory which is in its childhood, and this is the subject we shall explore during the rest of this article.


\section{Limit sets of discrete subgroups of $PSL(3,\C)$: an example }\label{s: example}

This section pretends to be a motivation for what is coming next. All concepts mentioned here are explained later with more care.

The group of automorphisms of $\P^n$ is the projective group $PSL(n+1,\C)$, and 
 from now on we are going to  focus our attention on its discrete subgroups. 
 
 Recall that in the classical case described above, the projective line $\P^1$ splits in two invariant sets: the discontinuity region and the limit set. In the first of these sets the action is properly discontinuous; this is also the region where the group forms a normal family, so it is the region of equicontinuity of the group. The quotient of the region of equicontinuity by the group is an orbifold with remarkable geometry, which has been the focus of study of many authors for decades.
 
 On the other hand,  the limit set is where the 
  dynamics concentrates. For nonelementary groups, each orbit in the limit set is dense in it, and every other orbit accumulates in it. There is a vast literature about the  study of the geometry and dynamics of the limit set.
  
  Most of these claims remain valid when we study conformal Kleinian groups in higher dimensions, and 
  one would like to make a similar study in higher dimensions, for the discrete subgroups of 
  $PSL(n+1,\C)$.  
  However, when trying to do so, one faces many new problems that do not appear in conformal geometry. The first of these being  that there is not a canonical 
concept of limit set  for the group,  as  the following example shows.

\medskip

Let  $\gamma$ be the projective transformation in $\P^2$  induced by the matrix:
\[
\left (
\begin{array}{lll}
2^{-1/2}&0 & 0 \\
0  &1& 0\\
 0& 0 &2
\end{array}
\right )
\]

We denote by $\Gamma$ the cyclic
subgroup of $PSL(3,\C)$  generated by $\gamma$. Notice that we have three linearly independent eigenvectors  $\{e_1,e_2, e_3\}$, and each of these determines a fixed point for the action in $\P^2$. For simplicity we denote the fixed points in $\P^2$ by the same symbols $ \{e_1,e_2, e_3\}$. Notice that $2^{-1/2} < 1 < 2$. Hence 
the point $\{e_1\}$ is  repellent,  a source, while $\{e_2\}$ is a saddle and $\{e_3\}$ is an attractor. The projective lines $\overleftrightarrow{e_1,e_2}$ and  $\overleftrightarrow{e_2,e_3}$ are both invariant lines. The orbits of points in the line $\overleftrightarrow{e_1,e_2}$  accumulate at $e_1$ going backwards, and they accumulate in $e_2$ going forwards. Similar considerations apply to the line $\overleftrightarrow{e_2,e_3}$, now $e_2$ is repellent in this line.

The
orbit of each point in $\P^2 \setminus (\overleftrightarrow{e_1,e_2} \cup \overleftrightarrow{e_2,e_3})$ 
accumulates at the points  $\{e_1,e_3\}$ and it is not hard to prove:

\begin{proposition} \label{e:inicia} If $\,\Gamma= \langle \gamma \rangle$ is the cyclic group generated by $\gamma$, then:
\begin{enumerate}
\item \label{f:e1} $\Gamma$ acts   discontinuously on the sets 
$\;\Omega_0=\P^2-(\overleftrightarrow{e_1,e_2}\cup
\overleftrightarrow{e_3,e_2})\;$,
$\,\Omega_1=\P^2-(\overleftrightarrow{e_1,e_2}\cup \{e_3\})\;$ and
$\,\Omega_2=\P^2-(\overleftrightarrow{e_3,e_2}\cup \{e_1\})$.

\item \label{f:e2} $\Omega_1$ and $\Omega_2$ are the maximal open
sets where $\Gamma$ acts (properly) discontinuously; and
$\Omega_1/\Gamma$ and $\Omega_2/\Gamma$ are compact complex
manifolds. (In fact they are Hopf manifolds).

\item \label{f:e3}$\Omega_0$  is the largest open set where
$\Gamma$ forms a normal family.
\end{enumerate}
\end{proposition}

 Thence, unlike the conformal case, now there is not a largest region where the action is properly discontinuous, but there are two maximal ones, $\,\Omega_1$ and $\,\Omega_2$, and none of them coincides with the 
  region of equicontinuity $\,\Omega_0$. 
Also,  
 now different orbits
accumulate at different points,  while  the set of accumulation points
of all orbits consists of the points $\{e_1, e_2,e_3\}$.  And 
it follows from the  proposition above that $\Gamma$ is not  acting
discontinuously on the complement of this set.

From this example we see that even in simple cases, when we
look at actions on higher dimensional projective spaces,  there is
not a definition of the limit set having all the properties one has in
the conformal setting. Here one might take as ``limit set":

\begin{itemize}
\item  The points $\{e_1,e_2, e_3\}$ where all orbits accumulate.
But the action  is not  properly discontinuous on all of its complement. Yet,
this definition is good in some sense. This corresponds
to taking the  Chen-Greenberg limit set of discrete subgroups of $PU(2,1) \subset PSL(3,\C)$.

\item The two lines $\overleftrightarrow{e_1,e_2}$ , $
\overleftrightarrow{e_3,e_2}$, which  are attractive sets for the
iterations of  $\gamma$ (in one case) or  $\gamma^{-1}$ (in the
other case). This corresponds to Kulkarni's limit set of $\Gamma$,
that we define below, and it has the nice property that the action
on its complement is discontinuous and also, in this case,
equicontinuous. And yet, the proposition above says that away from
either one of these two lines (and a point)  the action of $\Gamma $ is properly 
discontinuous. 

\item We may be tempted to taking as limit set the complement
of the ``largest  region where the action is properly discontinuous", but there is no such
region: there are two of them, $\,\Omega_1$ and $\,\Omega_2$, so which one we choose?

\item Similarly we may want to define the limit set as the
complement of ``the equicontinuity region" $\Omega_0$. In this particular
example, that definition may seem to be appropriate. Yet, there are interesting complex manifolds one gets as quotients of
$\Omega_1$ and $\Omega_2$ by the action of $\Gamma$, that  cannot be
obtained as quotient of the region of equicontinuity of some
subgroup of $PSL(3,\C)$. That is, these manifolds can not be written
in the form $U/G$ where $G$ is a discrete subgroup of $PSL(3,\C)$
acting equicontinuously on an open set $U$ of $\P^2$.
 Moreover,
 there are examples where $\Gamma$ is the fundamental group of certain
 Inoue surfaces  and
 the action of $\Gamma$   on $\P^2$ has no points of equicontinuity.
\end{itemize}

Thus one has different definitions of the concept of limit set, each
with its own  nice properties in different
settings.

So  we see that in the higher dimensional setting we cannot expect a well defined notion of limit set, in the sense that 
the equicontinuity set, maximal domains of discontinuity and ``wandering`` domains may  not  agree as in the one 
dimensional case. One way to overcome such a problem  is by considering all the ``interesting"  sets where the action is 
properly discontinuous  and study  the relation between them, 
its   geometry as well as their analytic properties.  This is what we do in the sequel.


\section{Complex Kleinian groups and discontinuity sets}

We recall that the complex projective space
$\mathbb{P}^n_{\mathbb{C}}$ is defined as:
$$ \mathbb{P}^{n}_{\mathbb{C}}=(\mathbb{C}^{n+1}- \{0\})/\sim \,,$$
where "$\sim$" denotes the equivalence relation given by $x\sim y$
if and only if \newline $x=\alpha y$ for some  nonzero complex scalar
$\alpha$. This is   a  compact connected  complex $n$-dimensional
manifold, diffeomorphic to the orbit space $\s^{2n+1}/U(1)$, where
$U(1)$ is acting coordinate-wise on the unit sphere in $\C^{n+1}$.

It is clear that every linear automorphism of $\C^{n+1}$ defines a
holomorphic automorphism of $\P^n$, and it is well-known  that every automorphism of $\P^n$ arises in
this way.  Thus one has:

\begin{theorem}
The group of automorphisms of $\P^n$
is the projective group:
$$PSL(n+1, \mathbb{C}) \,:=\, GL({n+1}, \C)/(\C^*)^{n+1} \cong SL({n+1}, \mathbb{C})/\mathbb{Z}_{n+1} \,,$$
where $(\C^*)^{n+1} $ is being regarded as the subgroup of diagonal
matrices with a single nonzero eigenvalue, and we consider the
action of $\mathbb{Z}_{n+1}$ (viewed as the roots of the unity) on
$SL(n+1, \mathbb{C})$ given by the usual scalar multiplication. 
\end{theorem}

We want to study discrete subgroups of $PSL(n+1, \mathbb{C})$ that generalize the notion of Kleinian groups described in the first section of this article. 
We know from the previous section that there is not a
``well-defined" concept of limit set for discrete subgroups of
$PSL(n+1,\C)$. And yet, there is a well-defined concept of
``discontinuous action":

\begin{definition}\label{Discontinuous action} Let $\Gamma$ be a discrete group  of $PSL(n+1,\Bbb{C})$ and
 consider its natural action on $\Bbb{P}^n_{\Bbb{C}}$.
Let  $\Omega$ be a
$\Gamma$-invariant open subset of $\P^n$. The action of $\Gamma$ is
{\it discontinuous} on $\Omega$ if each point $x \in \Omega$ has a
neighborhood $U_x$ which intersects at most finitely many copies of
its $\Gamma$-orbit. The discontinuity region of $\Gamma$ is the largest such set in $\P^n$, union of all open invariant sets where the action is discontinuous.
\end{definition} 

Alas, unlike the case of groups acting on $\P^1$, now the action on $\Omega$ may not be neither properly discontinuous  nor equicontinuous:

\begin{definition} 
Let  $\Omega$ be a
$\Gamma$-invariant open subset of $\P^n$. The action of $\Gamma$ is
{\it properly discontinuous} on $\Omega$ if  each compact set $K \subset \Omega$   intersects at most finitely many copies of
its $\Gamma$-orbit. 
\end{definition}

 For instance, consider the cyclic group of diffeomorphisms of $\R^2$ generated by $(x,y) \mapsto (2x, \frac{1}{2} y)$. Then the action is discontinuous on $\R^2 \setminus (0,0)$, but to get a properly discontinuous we must remove one of the axes. We have a similar picture in the example envisaged in the previous section.

\begin{definition}
Let $G$ be a group acting on a manifold $X$. {\it The equicontinuity
region} of $G$, in symbols $Eq(G)$, is defined to be the set of points
$z\in \mathbb{P}^n_\mathbb{C}$ for which there is an open
neighborhood $U$ of  $z$   such that $G \mid_U$  is a normal family.
\end{definition}

\medskip
Now we define the Kulkarni limit set. This definition has the advantage of granting that the action on the complement of this limit set is properly discontinuous.

For this, 
 let
 $L_0(\Gamma)$  be the
closure of the set of points in $\Bbb{P}^n_{\Bbb{C}}$ with infinite isotropy group.
Let $L_1(G)$ be the closure of the set of cluster points of orbits
of points in $X-L_0(G)$, {\it i.e.},  the cluster points of the
family $\{\gamma(x)\}_{\gamma\in{G}}$, where $x$ runs over
$X-L_0(G)$.
Finally, let $L_2(G)$ be the closure of the set of cluster points of
$\{\gamma(K)\}_{\gamma\in{G}}$, where $K$ runs over all the compact
subsets
of $X-\{L_0(G)\cup{L_1(G)}\}$. 
\medskip

For instance, it is an exercise to see that in the previous example, when $ \Gamma$ is the cyclic subgroup of $PSL(3,\C)$ determined by the matrix
\[
\left (
\begin{array}{lll}
2^{-1/2}&0 & 0 \\
0  &1& 0\\
 0& 0 &2
\end{array}
\right ) \;, 
\]
one has that $L_0(\Gamma)$ and ${L_1(\Gamma)}$ are equal and they consist of the three fixed points 
$\{e_1,e_2, e_3\}$, 
while the set $L_2(\Gamma)$ consists of the  invariant lines $\overleftrightarrow{e_1,e_2}$ and  $\overleftrightarrow{e_2,e_3}$.

\medskip

We have:

\begin{definition}
 Let $\Gamma$ be a group as above, we define  
the {\it Kulkarni limit set} of $G$ in $X$ as the set
 $$\Lambda_{Kul}(\Gamma):=L_0(\Gamma)\cup{L_1(\Gamma)}\cup{L_2(\Gamma)}.$$
We also define  the {\it Kulkarni region of
discontinuity} of $\Gamma$ as
 $$\Omega_{Kul}(\Gamma)\subset{X}:= X-\Lambda_{Kul}(\Gamma).$$
\end{definition}

Notice that the set $\Lambda_{Kul}(G)$ is closed and $G$-invariant, and 
  $G$ acts properly
discontinuously on its complement. 

We remark that this definition of a limit set given in \cite{Ku1} applies in a very general setting of discrete group actions. In  
the classical case, when the group is a subgroup of conformal diffeomorphisms of $\s^n$, one can prove that the {
discontinuity } and the equicontinuity sets  agree with the Kulkarni region of discontinuity.

\begin{definition}
A group $\Gamma\subset PSL(n+1,\Bbb{C})$ is complex Kleinian if there is an non-empty open set  $\Omega\subset \Bbb{P}^n_
\Bbb{C}$ on which  $\Gamma$ acts properly discontinuously.
\end{definition}

There are two special families of such groups: those coming from complex hyperbolic geometry, and those coming from complex affine geometry. In the first case, these are the projectivisation of discrete subgroups of $U(n,1)$, the group of linear transformations of $\C^{n+1}$ that preserve a quadratic form of signature $(n,1)$. In the second case one has groups of affine transformations of $\C^n$; these extend to projective automorphisms of $\P^n$ that preserve the hyperplane at infinity.

In the following sections we will describe  other types of complex Kleinian groups and some of their properties.


\section{Complex Kleinian groups in $PSL(3,\C)$} \label{sck2}

In this section we study discrete groups of  $PSL(3,\C)$ acting on $\P^2$ with a nonempty region of discontinuity.
We exhibit several types of examples of such groups,  which for the moment are the only known Kleinian  subgroups 
of $PSL(3,\C)$. Furthermore, one can show that  under 
certain conditions, there are no more subgroups of $PSL(3,\C)$ up to conjugation.

As mentioned above, 
the first type of complex Kleinian groups we consider are the discrete subgroups of $PU(2,1)$. This is the group 
of holomorphic isometries of the complex hyperbolic space $\mathbb H^2_\C$, and its discrete subgroups are 
known as  {\it complex hyperbolic groups}. 
  
  The complex hyperbolic space $\mathbb H^2_\C$ can be regarded as being the 4-ball $\mathbb B^4$ in $\P^2$ consisting of points whose homogeneous coordinates $[z_1:z_2:z_3]$ satisfy:
  $$ |z_1|^2 +  |z_2|^2 <  |z_3|^2 \,.$$
  Its boundary is a 3-sphere in $\P^2$, which is the projectivisation of the cone of light in $\C^3$ given by 
  $$  |z_1|^2 +  |z_2|^2 - |z_3|^2  = 0 \;.$$
  Just as in real hyperbolic geometry, this boundary sphere is known as the sphere at infinity, so we may denote it $\s^3_\infty$.

  This ball $\mathbb B^4 $ becomes the space $\mathbb H^2_\C$ when we equip it with the Bergmann metric, which is a Hermitian metric with constant negative holomorphic sectional curvature; $PU(2,1)$ becomes its group of holomorphic isometries.

   Given a discrete  group
 $\Gamma\subset PU(2,1)$, let  $\Lambda_{CG}(\Gamma)$ be the set of accumulation points of the orbits $\Gamma x$, for some $x\in \mathbb H^2_\C $. This is the {\it Chen-Greenberg} limit set of $\Gamma$. 
 
 Just as for classical Kleinian groups, since the action of $\Gamma$ on $\mathbb H^2_\C $ is by isometries, it is also properly discontinuous there. It follows that limit set $\Lambda_{CG}(\Gamma)$ is contained in the sphere at infinity
 $\s^3_\infty$.

  This limit set  has most of the nice properties of the limit set of classical Kleinian groups. For instance, if it has finite cardinality, then it consists of at most two points, and in this case the group is elementary.
  Otherwise, 
   it   is not hard to check that  $\Lambda_{CG}(\Gamma)$ does not depend on the choice of $x$, and it  is a closet $\Gamma$-invariant set which is a nowhere dense subset of $\partial (\mathbb H^2_\C)$.    
   
   \medskip
   
   The following example due to Dutenhefner-Gusevskii  shows how  intricate the Chern-Greenberg limit set can be. Here we just  sketch the construction of the group, see \cite{Gus} for more details:  
 
 \begin{example} 

Recall that a knot means (classically) a copy of the circle $\s^1$ embedded in a 3-sphere.
Consider a knot $K\subset \partial \Bbb{H}^2_{\Bbb{C}}$ and a finite collection $S = \{S_k,
S_k^\prime \}, k =1\ldots n ,$ of  spheres  placed along $K$, contained in $\partial{\Bbb{H}^2_{\Bbb{C}}}$, satisfying the following condition: there
is an enumeration $T_1, \ldots , T_{2n}$ of the spheres of this family
such that each $T_k$ lies outside all the others, except that $T_k$
and $T_{k+1}$ are tangent, for $k =1,\ldots, 2n-1$, and $T_{2n} $
and $T_1$ are tangent. Such a collection $S$ of spheres is called a {\it Heisenberg string of beads}. Let $g_k$ be elements from $PU(2, 1)$ such that
\begin{enumerate}
\item $g_k (S_k) = S_k^\prime$,

\item $g_k(Ext(S_k)) \subset Int(S_k^\prime)$,

\item $g_k$ maps the points of tangency of $S_k$ to the points of
tangency of $S_k^\prime$,
\end{enumerate}


It is proved in \cite{Gus} that one can choose Heisenberg string of beads and the generating maps $g_k$ so that 
the limit set of the group $\Gamma$ is a wild
knot. In fact the knots obtained in this way are not differentiable at any point, and they are self-similar.
\end{example}

The following result describes the relations amongst the different notions of limit set in the case of complex hyperbolic groups, see \cite{Nav1}.

\begin{theorem}
Let $\Gamma\subset PU(2,1)$ be a discrete group.  Then the limit set in Kulkarni' s sense  is the union of all projective lines tangent to the sphere $\s^3_\infty$ at a point in the Chen-Greenberg limit set. That is:
\[
\Lambda_{Kul}(\Gamma)=\bigcup_{z\in \Lambda_{CG}(\Gamma)}\ell_z
\]
where $\ell_z$ is the projective  line  tangent to $\partial (\mathbb H^2_\C)$ at $z$. Moreover, if  $\Lambda_{CG}(\Gamma)$ contains more than 2 points, then $\Omega_{Kul}(\Gamma)$ is the largest open set where $\Gamma$ acts properly discontinuously, and it is also the equicontinuity region of $\Gamma$ in $\P^2$.
\end{theorem}

 Let us consider now other types of examples. Most of these are affine groups (see \cite{cano2} for details).
 
\begin{example}
{\it Fundamental Groups of Complex Tori}.
Let    $g_i:\Bbb{C}^2\rightarrow \Bbb{C}^2$, $i=1...4$, be translations induced by $\Bbb{R}$-linearly independent 
vectors in $\Bbb{C}^2$ and  let $\Gamma= \langle g_1,\ldots,g_{4} \rangle$ be the group  generated by these maps. It is not hard to show that the action of $\Gamma$ 
can be extended to  $\Bbb{P}^2_\Bbb{C}$ and one has $ L_0(\Gamma)\cup L_1(\Gamma)=\Lambda(\Gamma),$
$Eq(\Gamma)=\Omega_{Kul}(\Gamma)=\Bbb{C}^2$. In this case  
$\Omega_{Kul}(\Gamma)$ is the largest open set on which  $\Gamma$ acts properly discontinuously.
\end{example}

\begin{example} {\it The  Suspensions}.
Given    $g\in PSL(2,\Bbb{C})$  we can induce  a map in $SL(3, \Bbb{C})$ as follows:
\[
\left (
\begin{array}{ll}
\tilde g & 0\\
0 & 1\\
\end{array}
\right);
\]
where $\tilde g$ is a lift of $g$. More generally, given   $\Sigma \subset \Bbb{C}^*$ a Kleinian group and $G\subset 
\Bbb{C}^*$ a  discrete group we define the suspension $Sus(\Sigma,G)$ of $\Sigma$ with respect to $G$ as
\[
\left \{
 \left (
\begin{array}{lll}
g \tilde h& 0\\
0 &  g^{-2}\\
\end{array}
\right )
 :g\in G, \tilde h \textrm{ is a lift of an element in } \Sigma \right\} \;.
 \]
Set
\[
\mathcal{C}= \left \{
\begin{array}{ll}
\bigcup_{p\in \Lambda(\Gamma)} \overleftrightarrow{p,e_3} & \textrm{if } G \textrm{ is finite.}\\
\bigcup_{p\in \Lambda(\Gamma)} \overleftrightarrow{p,e_3} \cup \overleftrightarrow{e_1,e_2} & \textrm{if } G
\textrm{ is infinite.}
\end{array}
\right.
\]
Then  $Eq(Sus(G, \Gamma))=\Bbb{P}^2_\Bbb{C}\setminus\mathcal{C}=\Omega_{Kul}(G)$  is the largest open 
set on which  $\Gamma$ acts properly discontinuously.
\end{example}

\begin{example} {\it Inoue's surfaces.}
Let  $M\in SL(3,\mathbb{Z})$ be a matrix with eigenvalues   $\alpha,\beta,\overline{\beta}$. Let $(a_1,a_2,a_3)$ 
be an eigenvector for  $\alpha$ and  $(b_1,b_2,b_3)$ an eigenvector for $\beta$. Set
\[
\gamma_0(w,z)=
\left
(
\begin{array}{lll}
\alpha & 0     & 0\\
0      & \beta & 0\\
0      & 0     & 1\\
\end{array}
\right );
\]
\[
\gamma_i(w,z)=
\left (
\begin{array}{lll}
1 & 0 & a_i\\
0 & 1 & b_i\\
0 & 0 & 1\\
\end{array} 
\right ) \;.
\]
It can be shown that the Kulkarni discontinuity region is
$$\Omega_{Kul}(G)=(\Bbb{H}^+\cup \Bbb{H}^-)\times
\mathbb{C} \;, $$
where $\Bbb{H}^+$ and $\Bbb{H}^-$ denote, respectively, the upper and the lower half planes in $\C$.
 This is now the largest open set on which $\Gamma$ acts properly discontinuously. The orbit space 
  $\Omega_{Kul}(G_M)/G_M$ consists of  two copies of  a 3-torus bundle over a circle and 
 $Eq(G_M)=\emptyset$.
 \end{example}

\begin{example} {\it A group Induced by a toral automorphism.}
Let   $M\in SL(2,\Bbb{Z})$ be given by
 $$M=\left (
\begin{array}{cc}
3 & 5\\
-5 & 8
\end{array}
\right ). $$
Consider the group $\Gamma_M$  induced by
\[
\begin{array}{l}
\left (
\begin{array}{cc}
M^k & b\\
0 & 1
\end{array}
\right ) \;,
\end{array}
\]
where $k\in \Bbb{Z}$ and $b\in \Bbb{Z}\times {Z}$. Then:
 $$\Omega_{Kul}(\Gamma_M)=Eq(\Gamma_M)=\bigcup_{i,j=0,1}(\Bbb{H}^{(-1)^i}\times \Bbb{H}^{(-1)^j})\,.$$  There are 
$p_1,p_2\in \Bbb{P}^2_\Bbb{C}$ and a ``circle" $\mathcal{C}\subset \Bbb{P}^2_{\Bbb{C}}$ such that the sets:
 \[
\Omega_i= \Bbb{P}^2_{\Bbb{C}} \setminus \bigcup_{p\in \mathcal{C}}\overleftrightarrow{p,p_i} \,,
 \]
are maximal open sets on which  $\Gamma_M$ acts properly discontinuously.
\end{example}

The following interesting example due to  R. Schwartz in \cite {schwartz1}  
 exhibits   surprinsinly rich  dynamics. This is generated  by means of the classical  Pappus' Theorem, which  we now recall:

\begin{theorem}[Pappus]
Let three points $p$, $b$, $q$ be contained in a straight line in the plane $\R^2$, and another three points $r$, $t$, $s$ be contained in another straight line. Then the three pairwise intersections  $\overleftrightarrow{bs}\cap \overleftrightarrow{tp}$, $\overleftrightarrow{ps}\cap\overleftrightarrow{rq}$, and $\overleftrightarrow{pt}\cap \overleftrightarrow{rb}$ are all contained in a  straight line.
\end{theorem}


Of course that Pappus' theorem can be regarded in the projective space $\Bbb{P}^2_\Bbb{R}$.

 It is proved in \cite{schwartz1} that this configuration of points and lines in $\mathbb P^2_\R$ gives rise  
to three projective transformations $\,\iota, \, \tau_1$ and $\tau_2$ in $PSL(3,\Bbb{R})$ satisfying  the following. Let $\mathcal L_1$ be the line containing the points   $p,b,q$, $\mathcal L_2$ is the line containing $r,s,t$ and $\mathcal L_3$ is the new line determined by Pappus theorem. 
The first transformation $\iota$ permutes $\mathcal L_1$ and $\mathcal L_2$.
The second transformation is denoted $\tau_1$, carries $\mathcal L_2$ into $\mathcal L_3$. The third   transformation  $\tau_2$ fixes $\mathcal L_2$ and it carries $\mathcal L_1$ into $\mathcal L_3$. Also, the transformations satisfy the following relations: $\iota^2=id,\tau_1\iota\tau_2=\iota, \tau_2\iota\tau_1=\iota, \tau_1\iota\tau_1=\tau_2, \tau_1\iota\tau_1=\tau_2, \tau_2\iota\tau_2=\tau_2$.

  \begin{theorem}
   The group $\Gamma_P$ generated by $\iota \, \tau_1$ and $\tau_2$ is isomorphic to the modular group $PSL(2,\Bbb{Z})$, and the $\Gamma_P$-orbit of the three lines $\mathcal L_1$, $\mathcal L_2$, $\mathcal L_3$ determines a curve in the dual projective space $\check {\mathbb P}^2_\R$, which generically is a fractal curve.
\end{theorem}


 Notice that the inclusion $PSL(3,\R) \hookrightarrow PSL(3,\C) $ allows us to think of the modular group $\Gamma_P$ constructed above as acting on $\P^2$. Then one has the following theorem from \cite{BCN2}:

\begin{theorem} The group $\Gamma_P$ is not conjugate in  $PSL(3,\C) $ to any affine or complex hyperbolic group, and its 
 Kulkarni limit set   in $\P^2$ is the union of the complexification of the real lines in ${\mathbb P}^2_\R$ appearing in the $\Gamma_P$-orbit of the lines $\mathcal L_1$, $\mathcal L_2$, $\mathcal L_3$. 
 \end{theorem}


The proof of the previous result depends strongly on the fact that we are iterating Pappu's theorem in the  real projective plane, so a natural question is:
\begin{question}
    Does  the result still holds if we replace the  lines and points appearing in Pappus's theorem by points and lines in $\P^2$ which cannot be taken in  a real projective plane? \end{question}

In other words, consider $p$, $b$, $q$  in a complex projective line in  $\P^2$, and another three points $r$, $t$, $s$ be contained in another complex projective line in  $\P^2$. Then one has a Pappus theorem in this setting: if we now consider the projective lines determined by the various points, then 
the three pairwise intersections  $\overleftrightarrow{bs}\cap \overleftrightarrow{tp}$, $\overleftrightarrow{ps}\cap\overleftrightarrow{rq}$, and $\overleftrightarrow{pt}\cap \overleftrightarrow{rb}$ are all contained in a complex projective  line. We may now mimic Schwartz' construction and obtain a subgroup of $PSL(3,\C)$ isomorphic to the modular group. The question is: what can we say about the Kulkarni limit set of this group?

\section{Quasi-cocompact Groups}
As mentioned before, one of the main reasons for studying discrete groups acting on complex projective spaces, is that if  $\Omega$ an open inv ariant set  where the action is properly discontinuous, then the orbit space $\Omega/\Gamma$ is a complex orbifold equipped with a projective structure. This is particularly interesting when the quotient 
$\Omega/\Gamma$ is compact.

One can prove (see the following section) that if the group is elementary, then its 
 limit set is either a line, a line and a point, 2 lines or 3 non concurrent lines. The group is said to be complex hyperbolic if it is conjugate to a subgroup of $PU(2,1)$, and it is  affine if it is conjugate to a subgroup of affine transformations in $\C^2$. 
 
 The following result from \cite{cano2} refines and extends a theorem proved by Kobayashi and Ochiai in \cite{kobayashi} for compact complex surfaces with a projective structure.

\begin{theorem} \label{t:main}
Let  $\Gamma\subset PSL(3,\mathbb{C})$ be a quasi-cocompact, then  $\Gamma$ is  virtually affine or complex 
hyperbolic.
\end{theorem}

In fact we can give a complete classification of the possible open sets in $\P^2$ that can appear as an invariant region for a group action as above:  

\begin{theorem} \label{t:main}
Let  $\Gamma\subset PSL(3,\mathbb{C})$ be a quasi-cocompact group which is not virtually cyclic. Then $\Omega_{Kul}(\Gamma)$  is the largest  open set on which     $\Gamma$ acts properly and discontinuously, and 
  $\Omega_{Kul}(\Gamma)$ is of   one the following types, up to biholomorphism:
  
  i)  $\mathbb{C}^2$; 
  
  ii) $\mathbb{C}\times \mathbb{C}^*$;
  
  iii)  $\mathbb{C}^*\times \mathbb{C}^*$;
  
  iv) $\mathbb{C} \times(\mathbb{H}^-\cup \mathbb{H}^+)$ , where $\mathbb{H}^-$ and  $\mathbb{H}^+$ are, respectively, the lower and the upper half spaces in $\C$;

v) $\mathbb D^2 \times \mathbb{C}^*$, where $\mathbb D^2$ is the open unit disc in $\C$; 

vi) The complex hyperbolic space $\Bbb{H}^2_\Bbb{C}$.
\end{theorem}

We also have the classification of the types of groups one has in each case, as well as of the corresponding orbifolds one gets as the space of orbits.  For simplicity, we only give below the description of the types of surfaces one gets.
We  refer to \cite{cano2} for more information, details and the proofs of these theorems.

\begin{theorem}
Let  $\Gamma\subset PSL(3,\Bbb{C})$   be  a quasi-cocompact group  and let $\Omega_{Kul}(\Gamma)$ be its Kulkarni region of discontinuity. Then, up to a finite (possibly ramified) covering, the orbit space 
 $\Omega_{Kul}(\Gamma)/\Gamma$  is  of the following types:
 
 i)  a complex torus;
 
 ii) Hopf surface;
 
 iii)  a Kodaira primary surface;
 
 iv)  one or two copies of 
an Inoue surface;

v)  a complex hyperbolic manifold; or 

vi)  a union (possibly countable) of elliptic affine surface, where at least one is compact.
\end{theorem}

This classification depends on the classification of the compact complex surfaces which admit a projective structure given in  \cite{kobayashi}, as well as the description of the corresponding  developing maps and their holonomy groups  (\cite{klingler1, klingler2}), and some considerations  about the corresponding Kulkarni's limit set.

\section{Discontinuity Regions; the 2-dimensional Case}
 
 We know that for classical Kleinian subgroups of $PSL(2,\C)$ one has that the limit set   either consists of at most two points, or else it has infinite cardinality. In this section we explore this type of considerations for complex Kleinian subgroups of $PSL(3,\C)$. We will see that in this case the complement to any region where the group acts properly discontinuously  must contain at least one projective line, and if the group does not have an affine subgroup with finite index, then the limit set in  Kulkarni's sense  has   infinitely many projective lines Furthermore,  in that case it coincides with the equicontinuity region and it  is the largest open set where the group acts properly discontinuously. We also got more specific results,  whose presentation requires the following terminology.
 
 \medskip
Throughout this section $\Gamma$ denotes a Kleinian subgroup of $PSL(3,\Bbb{C})$.

\begin{definition}  We denote by  $Dis(\Gamma)$ the collection of all open invariant sets $\Omega\subset \Bbb{P}^2$ on which  $\Gamma$ acts properly discontinuously and which are maximal in the sense that they are not contained in a larger invariant open set where the group acts 
 properly 
discontinuously.
\end{definition}

Notice that  by Zorn's lemma, the set $Dis(\Gamma)$  is non empty.

It is not hard to show that in the complement of each set where $\Gamma$ acts properly discontinuously there is at least one complex projective line (see \cite{BCN}).  It is useful to have   control over the number of lines  appearing in the complements of elements in $Dis(\Gamma)$, so  we introduce the following number:

\begin{definition} Define $Lin(\Gamma)$ to be the smallest number of projective lines in the complement of some element in $Dis(\Gamma)$. That is:
$$Lin(\Gamma)=min \big \{ Card \, \{\ell\in Gr(\Bbb{P}^2_{\Bbb{C}}) \, \big \vert \, \ell\subset \Bbb{P}^2\setminus \Omega\; ; \;  \Omega\in 
Dis(\Gamma) \}  \big \} \,,$$
where $Gr(\Bbb{P}^2_{\Bbb{C}})$ denotes the set of complex lines contained in $\Bbb{P}^2_{\Bbb{C}}$.
\end{definition}

For instance, if $\Gamma$ is the cyclic group in the example in Section \ref{s: example}, then $Dis(\Gamma)$ has exactly 2 elements and $Lin(\Gamma)=1$.
If $\Gamma$ is a suspension group as in Section \ref{sck2},   then $Dis(\Gamma)=\{\Omega_{Kul(\Gamma)}\}$ and $Lin(\Gamma)=1,2,3, \infty$, depending on  the groups $G$ and  $\Sigma$ which appear in the construction (see \cite{BCN} for details).

\medskip

On the other hand, in the classical theory of Kleinian groups in $\P^1$, as well as in the study of the dynamics of rational maps in one or more complex variables, the use of normal families and  equicontinuity play a key-role. In order to allow the use of this analytic viewpoint for  studying  discrete subgroups of $PSL(3,\C)$, it is useful
to get estimates on the  number of lines in general position contained  in the complement  of elements in $Dis(\Gamma)$:
We  say that a collection of lines $\mathcal{L}$  in $\P^2$ is {\it an array of lines in general position} if there are no three  lines in this collection which are concurrent.

\begin{definition}
For each  $\Omega \in Dis(\Gamma)$
we define $Lg( \Omega)$ to be the set  of all arrays of lines in general position contained in its complement. That is, if we let $Gr(\Bbb{P}^2_{\Bbb{C}} \setminus \Omega) $ denote the set of projective lines that do not meet $\Omega$, then:
$$Lg( \Omega)= \{\mathcal{L}\subset  Gr(\Bbb{P}^2_{\Bbb{C}} \setminus \Omega) \, \big \vert \, \mathcal{L} \; \textrm{is an array of lines in general position} \}.$$

\end{definition}

\medskip

\begin{definition}
Given  $\Omega \in Dis(\Gamma)$, 
define $Ling(\Omega)$ to be the maximal number of lines in an array of lines in general position  contained in the complement of $\Omega$.
We now define $Ling(\Gamma)$ to be the minimal number of lines in an array of lines in general position for all $\Omega \in Dis(\Gamma)$.
$$Ling(\Gamma)=min\{ Ling(\Omega) \, \vert \,
 \Omega\in Dis(\Gamma)\}.$$

\end{definition}

\medskip

For instance if $\Gamma$ is a suspension group as in Section \ref{sck2}, then $Ling(\Gamma)$ is either 2 or 3 depending on whether  the group $G$ which appears in its construction is, respectively, finite or infinite.

The theorem below is proved in  \cite{BCN}.

\begin{theorem}\label{Th: number of lines}
 The number     $Lin(\Gamma)$ is 1, 2, 3 or infinite.  Moreover, $Ling(\Gamma)$ is 
1, 2, 3, 4 or infinite.
\end{theorem}

\begin{theorem}\label{Th: enough lines}
Let  $\Gamma$ be  such that     $Lin(\Gamma)\geq 2$. Then there is a largest open set 
where $\Gamma$ acts properly discontinuously. Moreover, if  $Ling(\Gamma)\geq 3$, then  $\Omega_{Kul}(\Gamma)=\Omega(\Gamma)=Eq(\Gamma)$ 
and 
$$
\Lambda_{Kul}(\Gamma)=\bigcup_{g\in \Gamma} \Lambda_{Kul}(\langle g \rangle).
$$
\end{theorem}
\begin{question}
 If $Ling(\Gamma)=2$, is it truth that $\Omega_{Kul}(\Gamma)$ is the largest open set where $\Gamma$ acts properly 
discontinuously?
\end{question}

 It would also be interesting to answer the following two questions. For these let $\Gamma\subset PSL(3,\C)$ be a complex Kleinian group which is not  a finite extension of a cyclic group.
 
 \begin{question}
Is there  a non-empty set $\Omega(\Gamma)\subset \P^2$ which is the 
largest open set where $\Gamma$ acts properly discontinuously ?  
\end{question}

\begin{question}
 If the answer to the previous question is positive, how does the dimension Hausdorff of $\P^2\setminus \Omega(\Gamma)$ varies when we study analytic families of groups?
\end{question}

In the case when $\Gamma$ is   a finite extension of a cyclic group the answer to these questions follows from
 (see \cite{Nav2, CNS}).

\medskip

In general 
the 
computation of the numbers $Lin(\Gamma)$ and $Ling(\Gamma)$   is not an easy task.  The following result provides  a 
large variety of groups where these numbers are $\infty$.
 
\begin{theorem}
If  $\Gamma$ is not virtually affine, then  
$Ling(\Gamma)=\infty $.
 
\end{theorem}

Recall that virtually affine means that it contains a finite index subgroup which is affine.

The following sequel of results from  \cite{BCN} is an attempt to describe the geometry of the sets which appear as discontinuity sets for  complex Kleinian groups.

\begin{theorem} 
 If  $Ling(\Gamma)> 4$, then $\Omega_{Kul}(\Gamma)$ is a holomorphy 
domain. Moreover $\Lambda_{Kul}(\Gamma)$ is a union of complex lines.
\end{theorem}

We recall that an open set in a complex manifold is said to be {\it a holomorphy domain} if there is a holomorphic function on it that does not extend to a larger set.

\begin{definition}
A group $\Gamma\subset PSL(3,\C)$ is said to be a {\it Toral group} if there is    a hyperbolic toral 
automorphism  $A\in SL(2,\Bbb{Z})$   such that $\Gamma$ is conjugate to the group

\[
\Gamma_A=
\left \{
\left (
\begin{array}{lll}
A^k  & b \\
 0 & 1\\
\end{array}
\right ) \, \big \vert \, b\in M(1\times 2,\Bbb{Z}), \, k\in \Bbb{Z}  
\right \} \,.
\]
\end{definition}

We know from Theorem  \ref{Th: number of lines} that the number $Ling(\Gamma)$ is either 1, 2, 3, 4 or $\infty$, and Theorem \ref {Th: enough lines} says that if this number is at least 3, then the Kulkarni region of discontinuity coincides with the region of equicontinuity and it
is the largest set where the group acts properly discontinuously.  It is thus natural to try to classify the groups with few lines. One has:

\begin{theorem} \label{t:23}
Let   $\Gamma\subset PSL(3,\Bbb{C})$ be a Kleinian group.
\begin{enumerate}
\item \label{t231} If   $Ling(\Gamma)=2$,  then  $\Omega(\Gamma)=U\times \Bbb{C}^*$ where $U$ is either the equicontinuity region 
of some group in $PSL(2,\C)$ or $\P^1\setminus U$ is  a countable  family of lines in  $\C $ with the  following 
property: if $(\ell_n)$ is a sequence of distinct lines in   $\mathcal{C}$, then $\ell_n\cup \{\infty\}$ converges, in the 
Hausdorff metric sense,   to either $\infty$ or $\Bbb{R}\cup \infty$. 
\item  \label{t232} If  $Ling(\Gamma)=3$, then   $\Omega(\Gamma)=U\times \Bbb{C}^*$ where $U$ is the equicontinuity region 
of some group in $PSL(2,\C)$
 \item 
One has  $Ling(\Gamma)=4$  if and only if  $\Gamma$ has a toral subgroup with 
index at most  8. In this case   $\Omega_{Kul}(\Gamma) $  is given by: 
\[
\Omega_{Kul}(\Gamma)=\bigcup_{\epsilon_1,\epsilon_2=\pm 1}\Bbb{H}^{\epsilon_1}\times \Bbb{H}^{\epsilon_2}, 
\]
where $\Bbb{H}^{+1}$ and $\Bbb{H}^{-1}$ denote the half-upper and half-lower  semiplanes  respectively.
\end{enumerate}
\end{theorem}

\medskip

 These results give rise to several questions: 
 
\begin{question} What can we say when
    $Ling(\Gamma)=1$?, is  $\Lambda(\Gamma)$ necessarily  a complex line or a line and a point?. \end{question}

\begin{question}
Is it possible to  improve the conclusion of part (\ref{t231}) of Theorem \ref{t:23} by saying  ``where $U$ is the  
equicontinuity region 
of some group in $PSL(2,\C)$"?
\end{question}

\begin{question}
 Is there  an example of a Kleinian group without proper invariant subspaces  which is not  conjugate to a complex  
hyperbolic group or a group generated  by the Pappus's procedure?
\end{question}


\section{The Higher dimensional setting}
We now briefly envisage complex Kleinian subgroups of $PSL(n+1,\C)$ for $n > 2$. This is indeed a branch of mathematics which is in its childhood and very little is known.

There are two specially interesting subgroups of $PSL(n+1,\C)$, each having a rich theory of discrete subgroups. The first of these is the projective Lorentz group $PU(n,1)$. This is the group of projective automorphisms of $\P^n$ whose action preserves the complex $n$-ball $\mathbb B$ of points in $\P^n$ with homogeneous coordinates $[z_1:\cdots:z_{n+1}]$ satisfying $$\vert z_1 \vert^2 + \cdots + \vert z_n \vert^2 < \vert z_{n+1} \vert^2 \;.$$
This ball can be equipped with the Bergmann metric and so becomes a model for the complex hyperbolic space. Its group of holomorphic isometries is precisely $PU(n,1)$, and there is a rich theory about its discrete subgroups. Each of these is automatically complex Kleinian, since the action on the ball, being by isometries,  is necessarily properly 
properly discontinuous. The subgroups of $PU(n,1)$ are usually known as {\it complex hyperbolic groups}.

Another specially interesting subgroup of $PSL(n+1,\C)$ is the group Aff$(\C^n)$ of affine transformations of $\C^n$. recall that $\P^n$ can be obtained by compactifying  $\C^n$ by attaching to it a copy of $\P^{n-1}$ at $\infty$. Every affine transformation of $\C^n$ naturally extends to the hyperplane $\P^{n-1}$ at $\infty$ and so one has a canonical embedding Aff$(\C^n) \hookrightarrow PSL(n+1,\C)$. Thence every discrete subgroup of Aff$(\C^n)$ can be regarded as being a discrete subgroup of $PSL(n+1,\C)$.

Now, what about discrete subgroup of $PSL(n+1,\C)$ which are neither affine nor complex hyperbolic?  Very little is known, and the constructions in \cite{No, SV1, SV2, SV3} point in this direction, and so do several of M. Kato's articles listed in our bibliography. Let us now discuss briefly some of these constructions, and
for simplicity we focus on the case $n=3$, {\it i.e.}, automorphisms of $\P^3$. 

As in  the one dimensional setting,  one may now construct   Schottky 
groups of $PSL(4,\C)$. Let us begin  with their formal definition:

\begin{definition} \label{d:s}
A subgroup $\Gamma \subset   PSL(n+1,\C)$  is called a {\it Schottky
group} if:

\begin{enumerate}
\item There are  $2g$ , $g\geq 2$, opens sets   $R_1,\ldots,R_g$, $S_1,\ldots,S_g$  in
$\P^n$ with the  property  that:
\begin{enumerate}
\item each of these open sets  is the interior of  its closure; and
\item the closures  of the $2g$  open sets  are pair-wise disjoint.
\end{enumerate}
\vskip.1cm
\item  $\Gamma$  has a generating set $Gen(\Gamma)=\{\gamma_1,\ldots,\gamma_g\}$ such that for all $1 \leq j \leq 
g$ one has that:
  $$\gamma(R_j)=\P^n \setminus \overline{S_j}\,,$$  where the
  bar means topological closure.
\end{enumerate}
If in the previous definition  we allow that interiors of the sets  $R_1,\ldots,R_g$, $S_1,\ldots,S_g$ are pairwise disjoint but not necessarily  the closures, the  respective group is  called  a {\it Kissing-Schottky group}.
\end{definition}

As in the one dimensional case, Schottky  groups are free with $g$ generators and the  associated manifold is compact. The following result ensures the existence of such groups acting on  $\P^3$, and actually on all projective 
spaces with odd dimension:

\begin{theorem} [Seade-Nori-Verjovsky] See \cite{SV1, SV2, No} 
Given $n,g>2$ there is a Schottky group acting on $\Bbb{P}^{2n-1}$ with $g$ generators.
\end{theorem}

The  construction  of Schottky groups done  by Nori is   in the spirit of the classical construction of Schottky group in the one 
dimensional case, as explained above. The 
 construction  in  \cite{SV1, SV2}  is inspired by the geometric construction of conformal Schottky  groups in $Conf_+(\s^n)$ by means of inversions on spheres. Each $(n-1)$-sphere defines an inversion in $\s^n$ and the group generated by these inversions is discrete whenever  the ``generating" spheres are chosen to be  pair-wise disjoint. 
 
 So in this construction each $(n-1)$-sphere in $\s^n$ plays the role of a {\it mirror}, splitting the sphere $\s^n$ in two diffeomorphic halves which are interchanged by the involution, and these involutions generate the group. The idea in \cite{SV1, SV2} is the same, the difference being that the mirrors are now the boundaries of tubular neighbourhoods of copies of $\P^n$ in $\P^{2n+1}$.  Let us explain this for $n=1$.
 
 Consider a projective line $\ell \cong \s^2$ in $\P^3$ and a tubular neighborhood  $\mathcal N(\ell)$. 
  It is an exercise to show that every such line has trivial normal bundle and therefore the boundary of  $\mathcal N(\ell)$ is diffeomorphic to $\s^2 \times \s^3$, but this happens only in dimension $3$ and it is not important for what follows. What really matters is that  the boundary 
 $\partial \mathcal N(\ell)$ is a real hypersurface in $\P^3$ which splits the projective space in two diffeomorphic halves, the interior of $\mathcal N(\ell)$ and the complement of its closure.  
 Furthermore,  there exist automorphisms of $\P^3$ that interchange these two halves, the interior and the exterior of 
 $\mathcal N(\ell)$, leaving the boundary as an invariant set. 
 
 Now choose pair-wise disjoint projective lines $\ell_1, \cdots, \ell_g$ in $\P^3$, and sufficiently thin tubular neighborhoods of these lines, $\mathcal N(\ell_1), \cdots, \mathcal N(\ell_g)$. For each of these, choose an element
 $\gamma_i \in PSL(4,\C)$ interchanging the interior and the exterior of $\mathcal N(\ell_i)$. Then the group $\Gamma := \langle \gamma_1, \cdots, \gamma_g \rangle$ generated by these maps is a complex Kleinian subgroup of $PSL(4,\C)$.  Its index 2 subgroup consisting of words of even length acts freely on a region of discontinuity, with compact quotient. The orbit space is a compact, complex manifold with a very rich geometry: it is a Pretzel Twistor space in the sense of Penrose \cite{Pe3}.
 
 This construction generalizes to higher dimensions   by considering pair-wise disjoint copies of the projective space $\P^n$ in $\P^{2n+1}$. In \cite{SV3} the authors study the geometry, topology  and dynamics of these groups. It is proved that one has in this case an appropriate definition of the limit set $\Lambda(\Gamma)$ and the action on its complement $\Omega(\Gamma) := \P^{2n+1} \setminus \Lambda(\Gamma)$  is properly discontinuous. The limit set is of the form ${\mathcal C} \times \P^n$ and the action on it is ``transversely minimal" 
 in the sense that if we take the corresponding action on the Grassmanian of $n$-planes in $\P^{2n+1}$, then the limits set is the Cantor set ${\mathcal C} $ and the action on it is minimal, {\it i.e.} all orbits are dense. 
 
 In all these cases, the corresponding subgroup of words of even length acts freely on the open set $\Omega(\Gamma) $ and the orbits space $M_\Gamma$ is a compact complex manifold. The authors determine the topology of these manifolds and study their Kuranishi deformation space.

 This construction fails to produce examples in even dimensions, and actually 
 the main result in  \cite{cano1} ensures that  such kind of groups cannot  exist
in even dimensions, though one does have kissing-Schottky groups. In particular one has 
kissing-Schottky groups acting on $\P^2$ which are obtained by suspension of  (kissing)-Schottky group in $\P^1$.
In these cases one has a maximal region where the group acts properly discontinuously, and its complement has exactly two lines in general position:

\begin{question}
 Given  a kissing-Schottky  group $\Gamma\subset PSL(3,\Bbb{C})$,  is it true that $Ling(\Gamma)\leq 2$?. 
\end{question}

 In \cite{SV1, SV2} there is another construction of complex Kleinian groups that 
 uses the twistor  fibration over the spheres $\s^{2n}$. For $n=2$ this can be explained as follows. 
 
 Consider the two dimensional quaternionic space ${\mathcal H}^2$, which can be regarded as being $\C^4$ with an additional structure coming from quaternionic multiplication. One may consider the space of left quaternionic lines
   in 
 ${\mathcal H}^2$. Every such line $\mathcal L$ corresponds to a copy of ${\mathcal H} \cong \C^2$, so we may look at the set of complex lines in it; this forms a copy of $\P^1$. 
 
 Identifying each quaternionic line to a point  we obtain the 4-sphere $\s^4$,  
which can be thought of as being the quaternionic projective line $\mathbb P^1_{\mathcal H}$. 
  Notice that this defines a projection
  $\P^3 \to \s^4$ where the fibre over each point is the set of complex lines in the given quaternionic line. So the fibre is $\P^1$. This is a fibre bundle  known as the Calabi-Penrose fibration, and also as the {\it  twistor fibration}. The fibres are known as the {\it twistor lines}.
  
  It is well known that the conformal group $Conf_+(\s^4)$ embeds canonically in the projective group $PSL(4,\C)$. This implies that if $\Gamma$ is a discrete group of orientation preserving conformal diffeomorphisms of the 4-sphere, then $\Gamma$  is canonically a group of automorphisms of $\P^3$, and it acts on $\P^3$ carrying twistor lines into twistor lines.
   
  In \cite{SV2} the authors prove that this action on $\P^3$ is by isometries when restricted to each twistor line, which is equipped with the usual round metric.  One gets as corollary that the Kulkarni limit set of the action on $\P^3$ is the lifting of the limit set of the action on $\s^4$. Hence the  dynamics of   discrete subgroups of conformal automorphisms of $\s^4$ embeds in the holomorphic dynamics of subgroups of $PSL(4,\C)$.

For instance, this implies that the fundamental group of every closed hyperbolic 5-manifold acts canonically on $\P^3$ and every orbit is dense. 

In particular, given two Kleinian groups acting on $\s^4$ we can use Maskit's combination theorems to produce a new Kleinian group in $Conf_+(\s^4)$, which lifts to a Kleinian group in $PSL(3,\C)$.  Recall that Maskit combination theorems  are one of the most important means for producing rich dynamics in the classical setting.
 The following  problem would produce very interesting  holomorphic dynamics in higher dimensions.
 
 \begin{question}
 Is it possible to give a generalization of the Maskit combination theorems to the case of complex Kleinian groups?
 \end{question}

 There are several interesting papers by M. Kato for subgroups of $\P^3$, including results in the spirit of Maskit theorems (see the references below).

\bibliographystyle{amsplain}

\begin{thebibliography}{10}

\bibitem{fatou}
P. Appell, E. \'I. Goursat,
{\it Th\'eorie des fonctions algebriques et de leurs integrales:
\'etude des fonctions analytiques sur une surface de Riemann},
1895,  Gauthier-Villars 


\bibitem{BCN}  W. Barrera, A. Cano, J. P. Navarrete. {\it Estimates of the number of lines lying in the limit set for subgroups of $PSL(3,\Bbb{C})$}, preprint, 2010, arXiv:1003.0708. 

\bibitem{BCN2}  W. Barrera, A. Cano, J. P. Navarrete. {\it A 2-dimensional complex Kleinian group with infinite lines in the limit set lying in general position}.
Preprint, 2010, arXiv:1003.0380.


\bibitem{ball}
W. W. R. Ball, H. S.  Coxeter, {\it Mathematical Recreations and Essays}, 13th ed. New York: Dover, pp. 155-161, 
1987. 

\bibitem{crash}
L. Bers, I. Kra, {\it A crash course on Kleinian groups}, Lecture Notes in Mathematics,  Springer-Verlag, 1974.




\bibitem{cano1}
A. Cano. {\em Schottky Groups are not Realizable in
$PSL(2n+1,\mathbb{C})$}. preprint, 2008.

\bibitem{cano2}
A. Cano. {\em Discrete Subgroups of automorphism of
$\mathbb{P}^2_{\mathbb{C}}$}. preprint, 2008.

\bibitem{CNS} A. Cano, J. P. Navarrete, J. Seade. {\em Complex Kleinian Groups}. Monography 2011, to be published.

\bibitem{Ca-Se} A. Cano, J. Seade. {\em On the equicontinuity region of discrete subgroups of $PU(n,1)$}. Preprint 
2008.

\bibitem{ChG}
S. S. Chen and L. Greenberg. {\em Hyperbolic Spaces}. Contributions
to Analysis, Academic Press, New York, pp. 49-87 (1974).

\bibitem{cox1}
H. S. M.Coxeter, {\it The Partition of a Sphere According to the Icosahedral Group}, Scripta Math 4, 156-157, 1936. 

\bibitem{cox2}
H. S. M.Coxeter,
{\it Regular Polytopes}, 3rd ed. New York: Dover, 1973. 


\bibitem{Gus}
F. Dutenhefner, N. Gusevskii, {\em Complex hyperbolic Kleinian
groups with limit set a wild knot}.  Topology  {\bf 43}  (2004),
no. 3, 677-696.



\bibitem{gray}
J. J. Gray, {\it Linear differential equations and group theory from Riemann to Poincar\'e`}, Birkh\"auser, 2000. 

\bibitem{kapovich}
M. Kapovich, {\it Kleinian Groups in Higher Dimensions}, Progress in Mathematics, Vol. 265, 485-562
2007. Birkh\"auser Verlag Basel/Switzerland.

\bibitem{Ka1}   M. Kato,    {\em On compact complex 3-folds with
lines}. Japanese J. Math.  {\bf  11} (1985),   1-58.

\bibitem{Ka2}   M. Kato,    {\em Factorization of compact
complex 3-folds which admit certain projective structures}. T\^ohoku
Math. J.  {\bf 41} (1989),   359-397.

\bibitem{Ka3}  M. Kato,   {\em Compact Complex 3-folds with
Projective Structures; The infinite Cyclic Fundamental Group Case}.
Saitama Math. J.  {\bf  4} (1986),  35-49.

\bibitem{Ka4}   M. Kato,    {\em Compact Quotient Manifolds of
Domains in a Complex 3-Dimensional Projective Space and the Lebesgue
Measure of Limit Sets}.   Tokyo J. Math.  {\bf  19} (1996), 99-119.

\bibitem{Ka5}   M. Kato, {\it
Compact quotients with positive algebraic dimensions
of large domains in a complex projective 3-space}, The Mathematical Society of Japan
J. Math. Soc. Japan
Vol. 62, No. 4 (2010) pp. 1317-1371.

\bibitem{kobayashi}  S. Kobayashi, T. Ochiai. {\em Holomorphic
Projective Structures on Compact Complex Surfaces (I and II)}. Math.
Ann.  {\bf 249}  (255), 1980 (81),
 75-94 (519-521).

\bibitem{klingler1}
B. Klingler. {\em Structures  Affines et Projectives sur les
Surfaces Complexes}. Annales de L' Institut Fourier, Grenoble, {\bf
48} 2 (1998),  441-447.


\bibitem{klingler2}
B. Klingler. {\em Un th\'eoreme de rigidit\'e non-m\'etrique pour
les variet\'es localement sym\'etriques hermitiennes}.   Comment.
Math. Helv.  {\bf  76}  (2001),  no. 2, 200-217.

\bibitem{krai}
M. Kraitchik,{\it  A Mosaic on the Sphere}, Mathematical Recreations, New York: W. W. Norton, pp. 208-209, 1942

\bibitem{Ku1}   R. S. Kulkarni.    {\em Groups with domains of
discontinuity}.  Math. Ann. {\bf 237} (1978),  253-272.

\bibitem{levi}
Ravi P. Agarwal, Donal O'Regan, {\it
 Ordinary and partial differential equations: with special functions, Fourier series, and boundary value problems}, Springer UGT, 2009.


\bibitem{maskit}
B. Maskit, {\it Kleinian groups}, Springer-Verlag, 1988

\bibitem{No}    M. V. Nori.    {\em The Schottky groups in higher
dimensions}.    Proceedings of Lefschetz Centennial Conference,
Mexico City, AMS Contemporary  Maths.{\bf  58}, part I (1986),
195-197.

\bibitem{indra}
 D. Mumford, C. Series and D. Wright,
{\it Indra's Pearls: The Vision of Felix Klein is a geometry}, y Cambridge University Press,  2002.

\bibitem{Nav1}   J.-P. Navarrete.  {\em
On the limit set of discrete subgroups of $\text{PU}(2,1)$}. Geom.
Dedicata 122, 1-13 (2006).

\bibitem{Nav2}   J.-P. Navarrete.  {\em   The Trace Function and Complex  Kleinian Groups in $\mathbb{P}^2_{\mathbb{C}}$}. Int. Jour. of
Math. Vol. 19, No. 7 (2008) 865-890.

\bibitem {Pe3}   R. Penrose.    {\em Pretzel twistor spaces}.
In ``Further advances in twistor theory'', ed. L.J. Mason and L. P.
Hughston, Pitman Research Notes in Maths.  vol. {\bf 231} (1990),
246-253.

\bibitem{poincare}
H. Poincar\'e, {\it M\'emoire sur Les groupes klein\'eens}, Acta Mathematica, 1883.

\bibitem{schwartz1}
R. Schwartz, {\it  Pappus' theorem and the modular group},  Inst. Hautes \'Etudes Sci. Publ. Math.  No. 78  (1993), pp. 187-206.


\bibitem{Seade}
J. Seade. {\em 
On the topology of isolated singularities in analytic spaces}. Progress in Mathematics 241. Birkh\"auser 2006.

\bibitem{SV1}
J. Seade, A. Verjovsky. {\em Actions of Discrete Groups  on Complex
Projective Spaces}. Contemporary Math.,  {\bf 269} (2001),155-178.

\bibitem{SV2}
J. Seade, A. Verjovsky. {\em Higher Dimensional Kleinian groups}.
Math. Ann., {\bf 322}  (2002), 279-300.

\bibitem{SV3}
J. Seade, A. Verjovsky. {\em Complex Schottky groups}. Geometric
methods in dynamics. II.  Ast\'erisque {\bf 287}  (2003), xx,
251-272.

\bibitem{map}
J. P. Snyder,  {\it Map Projections - A Working Manual}, ,1987. 
 
 \bibitem{Sul} D.  Sullivan. {em 
Quasiconformal homeomorphisms and dynamics. I: Solution of the Fatou- Julia problem on wandering domains}. Ann. Math. (2) {\bf 122} (1985), 401-418.
 
 \bibitem{thurston}
W. P. Thurston, {\it Three-Dimensional Geometry and Topology:
Volume 1}, Princeton Mathematical Series
, 2004.

\bibitem{Yo}
M. Yoshida, {\it Orbifold-uniformizing differential equations},
Mathematische Annalen
Volume 267, Number 1, 125-142.

\bibitem{love}
M. Yoshida, {\it Hypergeometric Functions, My Love: Modular Interpretations of Configuration Spaces}, Braunschweig-
Wiesbaden: Friedr. Vieweg \& Sohn, 1997.



\bibitem{galli}
J. Wyss-Gallifent, {\it Discreteness and Indiscreteness Results for Complex Hyperbolic
Triangle Groups}, Ph.D. Thesis, University of Maryland (2000).
\end{thebibliography}

\end{document}